\documentclass[letterpaper, 10 pt, conference]{ieeeconf}  

\IEEEoverridecommandlockouts                              

\title{\LARGE \bf
Error Bounds in Nonlinear Model Predictive Control with Linear Differential Inclusions of Parametric-Varying Embeddings 
}

\author{Dimitrios S. Karachalios$^{1}$, Maryam Nezami$^{1}$, Georg Schildbach$^{1}$ and Hossameldin S. Abbas$^{1}$
\thanks{*First author was supported by DFG (German research foundation)}
\thanks{$^{1}$The authors are with the Faculty of Electrical Engineering in Medicine, University of Luebeck, Germany. {\tt\small email:~dimitrios.karachalios@uni-luebeck.de}}%
}

\newtheorem{theorem}{Theorem}[section]

\newtheorem{remark}{Remark}[section]

\newtheorem{example}{Example}[section]
\def\IR{{\mathbb R}}

\def\IZ{{\mathbb Z}}

\newcommand{\cP}{ {\cal P} }

\def\IR{{\mathbb R}}

\usepackage{graphicx}
\usepackage{amsmath}
\usepackage{amsfonts}
\usepackage{amssymb}
\usepackage{svg}
\usepackage{times} 
\usepackage[capitalise]{cleveref}
\usepackage{algpseudocode}
\usepackage{algorithm}
\usepackage{caption}

\DeclareMathOperator{\sinc}{sinc}
\newtheorem{problem}{Problem}
\begin{document}
\maketitle
\thispagestyle{empty}
\pagestyle{empty}
\begin{abstract}
In this work, we provide deterministic error bounds for the actual state evolution of nonlinear systems embedded with the linear parametric variable (LPV) formulation and steered by model predictive control (MPC). The main novelty concerns the explicit derivation of these deterministic bounds as polytopic tubes using linear differential inclusions (LDIs), which provide exact error formulations compared to linearization schemes that introduce additional error and deteriorate conservatism. The analysis and method are certified by solving the regulator problem of an unbalanced disk that stands as a classical control benchmark example.
\end{abstract}
\section{INTRODUCTION}
Interesting control applications involve dynamics that do not satisfy the superposition and scaling principles thus, are nonlinear. On one hand, modeling with partial differential equations and discretizing with finite methods results in large-scale ordinary differential equations that only efficient model reduction can lead to feasible control \cite{AntoulasBook2005}. On the other hand, good observables result in highly nonlinear models of low-dimension easier to handle. Both approaches introduce errors due to discretization, reduction schemes, or model mismatches due to assumptions with the original plant. Handling wisely the error that comes from the modeling mismatches between the predictive model and the actual plant's response has driven research to provide robust-tube control methods for linear and nonlinear systems in a series of papers \cite{Mayne2011tubes,MAYNE2005219,Cannon2011,LANGSON2004125,Schildbach2013}. These methods perform well under disturbances or measurement noise. 

Among many nonlinear control methods, increasing attention has been drawn to nonlinear model predictive control (NMPC), which can cast the control task into an optimization problem that can handle in addition input-state constraints. Representing control as an optimization problem results in finding the global minimum over a nonconvex manifold that is quite challenging to solve in real-time. In addition, the nonlinear constraints in both input and states increase the complexity of the admissible search space for detecting that minimum. To overcome the aforementioned problems in NMPC, an excellent alternative is to embed the nonlinear system in a linear parameter-varying (LPV) formulation that results in a quadratic manifold with an inherent unique optimal solution that can be solved efficiently and online. Moreover, the nonlinear constraints can be realized with tangential schemes that allow adaptive linear constraint formulation. Combining the quadratic manifold along with the adaptive linear constraints, the linear parametric-varying model predictive control (LPVMPC) problem results in solving a classical quadratic program (QP) with many efficient algorithms to handle this operation in real-time as in autonomous driving tasks~\cite{Nezami2023}.

\begin{figure}
    \centering
    \includegraphics[width=0.55\textwidth]{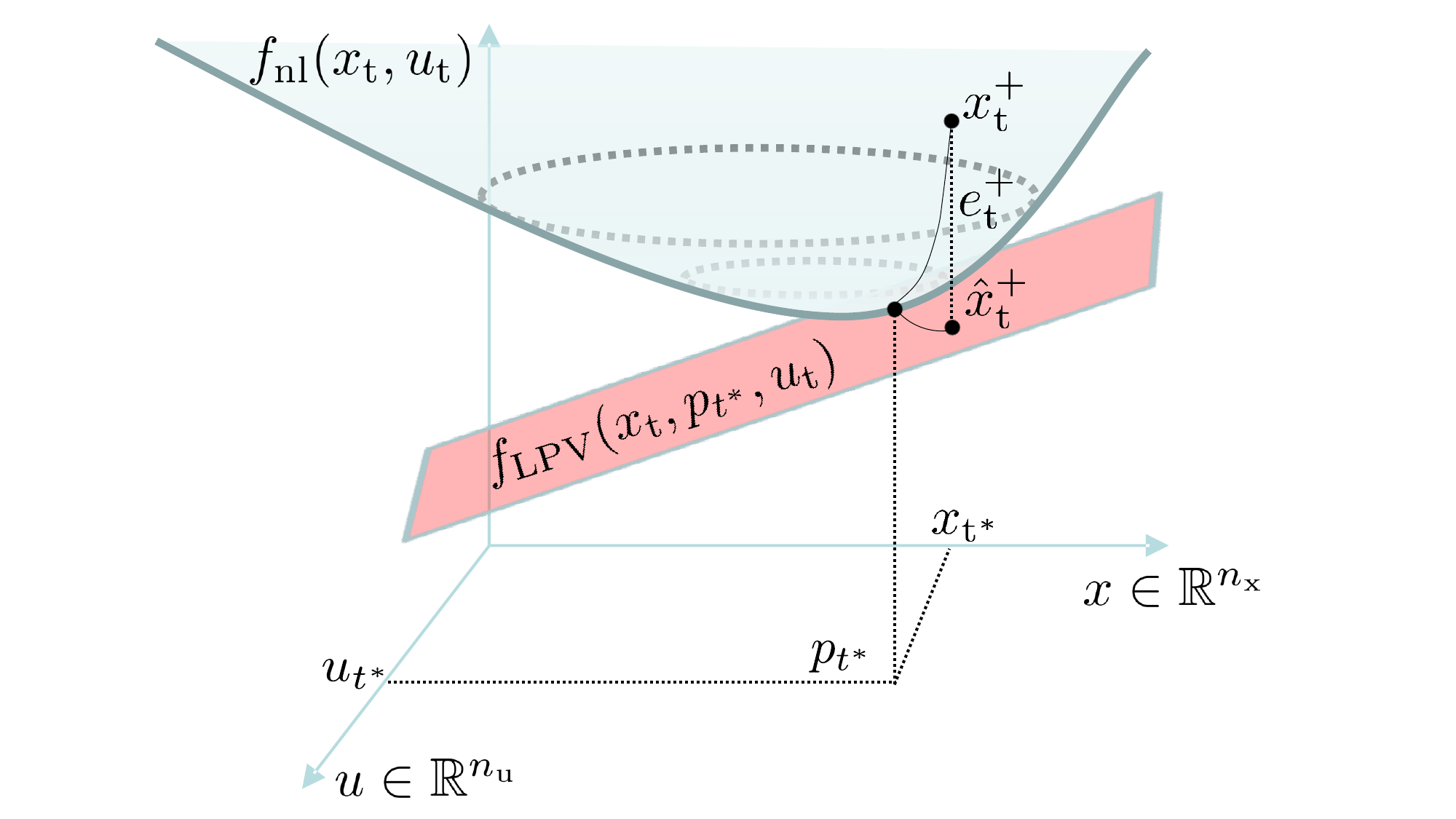}
    \caption{The LPV operator $f_{\text{LPV}}(x_t,u_t):=A(p_t^*)x_t+Bu_t=x_t^+$ realizes the nonlinear manifold with a sliding tangent hyperplane over the time $t$ under a fixed scheduling parameter $p_t^*$ that introduces the error state $e_t^+$.}
    \label{fig:LPV}
    \vspace{-6mm}
\end{figure}

The challenge in representing nonlinear systems with the LPV embedding is the appropriate prediction of the uncertainty quantification introduced by the so-called scheduling parameter $p$. The vector $p$ absorbs all the nonlinear dependencies in an affine representation, and the LPV operator realizes the nonlinear manifold with an adaptive tangential hyperplane for fixed values of $p$ \cref{fig:LPV}. Therefore, to solve the LPVMPC problem with QP, the prediction of the scheduling parameter within the receding horizon needs to be initialized prior and consequently produces an error between the actual and model responses.
The problem of the uncertainty related to embedding LPV models into an MPC framework has recently received considerable attention. For instance, one approach is the adoption of tube-based MPC for the LPV setting, as proposed in~\cite{HanemaTubes}. In this paper, under the assumption of the boundedness of future $p$, an anticipative tube MPC algorithm for LPV systems is proposed. Also, in~\cite{nezami2022robust}, an approach inspired by the anticipative tube MPC introduced in~\cite{HanemaTubes} for autonomous lane keeping is presented. This paper suggests a cascade MPC control architecture incorporating both LPV lateral and linear longitudinal models. Furthermore, in~\cite{Abbas2016}, the stability of the LPVMPC design is guaranteed by using a bilinear matrix inequality. However, this method can be conservative and computationally demanding. In~\cite{bujarbaruah2022robust}, an offline approach is proposed for finding bounds on model uncertainty as a function of the control inputs. In other words, the tightenings of the constraints in the MPC are functions of decision variables, i.e., control inputs. In addition to MPC-related studies, there are also some works that address the problem of providing bounds on LPV modeling mismatches. In~\cite{BAINIER202267}, the paper norm bounds the state trajectory of a continuous time LPV system. However, the nonlinear system needs to satisfy the Lipschitz condition. In~\cite{Chaillou2007}, by using the context of strongly monotone operators, an upper bound and a lower bound for a specific class of nonlinear systems is suggested. 

In this paper, at first, the derivation of the dynamics of the error for a general LPV system is explained. In the next step, we provide exact error bounds, between the LPV model and the actual nonlinear system, by utilizing the linear differential inclusion (LDI), \cite{BoydLDI}, of the nonlinear operator over the parametric uncertainty that is produced in the LPVMPC framework. The proposed bound on error can serve as a foundation for establishing further theoretical guarantees, e.g., stability, and recursive feasibility, which will justify safety features in control of the system. 


In \cref{sec:Prel}, we start with preliminaries and formal representation of the problem under consideration. In \cref{sec:Analdyna}, we formulate explicitly the error dynamics that evolve within the receding horizon along with the polytopic error bounds. In \cref{sec:results}, we apply our method to a classical control benchmark with ease of reproducibility. Finally, in \cref{sec:conclusion}, we summarize our findings and provide the open challenges and future research directions.
\section{Preliminaries \& Problem formulation}\label{sec:Prel}
\subsection{Definitions \& assumptions}
We start with the discrete nonlinear dynamical system 
\begin{equation}\label{eq:sysnlc}
    \Sigma:x_{k+1}=f(x_k,u_k),
\end{equation}
of state dimension $n_{\mathrm{x}}$, and input dimension $n_u$. Considering the sampling time $t_s$, it holds $t_k=t_sk,~\forall k\in\IZ_+$, with $x_k=x(t_sk)$, and $x_0=x(0)$ the initial condition state vector. $f:\IR^{n_{\mathrm{x}}}\times\IR^{n_u}\rightarrow\IR^{n_{\mathrm{x}}}$ is a nonlinear operator. The nonlinear dynamical system in \eqref{eq:sysnlc} can be represented equivalently with a parametric varying linear (LPV) formulation that will give rise to methods that do linear in an adaptive way. An appropriate scheduling parameter vector $p$ with dimension $n_{p}$ should be introduced to accomplish that. Towards one further simplification, the remaining LPV system and through a filter that will increase the state dimension can recast the scheduling dependence only to the linear matrix $A(p)$ allowing a static $B$. As a result, the original nonlinear system \eqref{eq:sysnlc} can be embedded in the following linear parameter varying (LPV) formulation \eqref{sys:LPV} as
\begin{equation}\label{sys:LPV}
    \Sigma:\left\{\begin{aligned}
        x_{k+1}&=A(p_k)x_k+Bu_k,\\
              p_k&=\rho(x_k,u_k),~x_0=x(0),\end{aligned}\right.
\end{equation}
where the mapping $\rho:\IR^{n_{\mathrm{x}}}\times\IR^{n_u}\rightarrow\IR^{n_{p}}$ is also given explicitly. In particular, the $\rho(\cdot)$ is a known nonlinear function of the state/input-$(x,u)$, which allows the embedding of \eqref{eq:sysnlc} in \eqref{sys:LPV}. Furthermore, the following parameterized matrix $A(p_k):\IR^{n_p}\rightarrow\IR^{n\times n}$ is also known and affine in terms of the scheduling parameter $p$. In particular, the affine structure of the discrete operator $A(p_k)$ can be expressed as
\begin{equation}\label{eq:affine}
    A(p_k):=A_0+\sum_{l=1}^{n_p}p_k^{[l]}A_{l},
\end{equation}
where $p_k^{[l]}$ denotes the $l^{th}$-element of the vector $p$ and $A_l$ are constant matrices. Together with the input matrix $B\in\IR^{n_{\mathrm{x}}\times n_u}$, the discrete-time LPV system is well-defined. The following remark \ref{rem:stanassu} summarizes the appropriate general assumptions to proceed with what follows for solving the control task by considering LPV predictive models.
\begin{remark}[Standing assumptions]\label{rem:stanassu}
To reflect the generalization of the method, here are the minimal assumptions:
    \begin{itemize}
        \item Appropriately smoothness of the nonlinear operator $f$ has been assumed (i.e., higher-order differentiability and continuity).
        \item The scheduling parameter can be measured at each sampling time $k$ but remains unknown within the receding horizon of length N.
        \item The input matrix $B$ does not depend on the scheduling parameter. This can be relaxed easily through a p-filter.
        \item No other disturbances or measurement noise has been assumed to affect the system $\Sigma$.
    \end{itemize}
\end{remark}
\begin{problem}[Error propagation within the receding horizon]
    We are interested in bounding the error produced between the true response of the actual nonlinear system and that predicted via LPV within the receding horizon control strategy (i.e., MPC) under a fixed scheduling prediction signal \cref{fig:LPV} that is inferred from prior knowledge.
\end{problem}
\subsection{Stabilization \& model predictive control} 
Unstable operation modes characterize the underlying dynamics in many interesting control applications. Therefore, control strategies that stabilize and drive the system to desired states under constraints are crucial. In the linear case, linear quadratic regulation (LQR) as a feedback state (i.e., $u=Kx$) efficiently provides stable closed-loop systems. In addition, LQR has been extended to handle LPV representations and provides robust-model state feedback \cite{robustLPVLQR} that stabilizes the system for all possible parametrization of the scheduling $p\in\cP$—a drawback, though, for LQR is that it cannot handle input-state constraints. Thus, we propose splitting the control input into two parts; the first part will stabilize the nonlinear plant through an LPV-LQR controller, and the second part will handle the input-state constraints. Specifically, the input design has the following structure:
\begin{equation}\label{eq:LQR}
    u_k=u_k^{\text{LQR}}+u_k^{\text{MPC}}=K x_k+u_k^{\text{MPC}}.
\end{equation}
By substituting \eqref{eq:LQR} to \eqref{sys:LPV}, we result to
\begin{equation}
\footnotesize
    \begin{aligned}
        x_{k+1}&=A(p_k)x_k+Bu_k,\\
        &=A(p_k)x_k+B(K x_k+u_k^{\text{MPC}}),\\
         &=\underbrace{\left(A(p_k)+BK\right)}_{A_c(p_k)}x_k+Bu_k^{\text{MPC}}.\\
    \end{aligned}
\end{equation}
Thus, the closed-loop dynamics of applying the LQR controller makes $A_c(p_k)$ stable and can be written as
\begin{equation}\label{sys:LPVLQR}
    \Sigma_{\text{c}}:\left\{\begin{aligned}
        x_{k+1}&=A_c(p_k)x_k+Bu_k,\\
              p_k&=\rho(x_k,u_k),
    \end{aligned}\right.
\end{equation}
with $A_c(p_k):=A(p_k)+BK$, where we denote with subscript ``c" the closed-loop operator and the remaining controller $u_k ^{\text{MPC}}$ can be denoted again as $u_k$ without asserting any confusion.
\subsection{Model predictive control (MPC) with LPV embeddings} After stabilizing the dynamics with the LQR feedback control, we want to drive the system to a given reference under some input and state constraints. Thus, we use MPC for control tasks with input or state constraints. The whole control problem can be cast as a constrained optimization problem within a given receding horizon of length $N$. The energy (cost) can be penalized with the quadratic weighted matrices\footnote{The quadratic weighted cost is defined as $\lVert x\rVert_Q^2=x^{\top}Qx$. Similarly, for $R$ and $P$.} $Q,~R$ and the quadratic cost $P$ for the terminal cost. At $t_k=k\cdot t_s$ and for $i=0,\ldots,N-1$ that will lead to a classical quadratic program (QP) along with the efficient algorithms that will allow real-time performance. The energy function to be minimized is
\begin{equation}
\footnotesize
\begin{aligned}
    J_k(u_{i|k})&:=\underset{u_{i|k}^*}{\min}\sum_{i=0}^{N-1}\left(\lVert \hat{x}_{i|k}\rVert_Q^2+\lVert u_{i|k}\rVert_R^2\right)+\lVert \hat{x}_{N|k}\rVert_P^2.
\end{aligned}
\end{equation}
In addition, we can introduce adaptive linear constraints. These are input and state constraints and can be introduced with the following sets:
\begin{equation}\label{eq:constraints}
    \begin{aligned}
        &\hat{x}_{i|k}\in\mathcal{X}_{i|k}=\{\hat{x}_k\in\IR^n|G_k^x \hat{x}_k\leq h_k^x\},\\
        &u_{i|k}\in\mathcal{U}_{i|k}=\{u_k\in\IR^{m}|G_k^u u_k\leq h_k^u\}.\\
    \end{aligned}
\end{equation}
\begin{problem}{QP optimization as $\texttt{QP}(\hat{p}_{i|k},\hat{x}_k,x_{i|k}^{\text{ref}})$}\label{prob:Prob_LPVMPC}
\begin{subequations}\label{eq:LPV_MPC}
\footnotesize
	\begin{align} 
		\underset{u_{i|k}^*}{\text{min}} \
		& \! \sum_{i=0}^{N-1}\left(\lVert \hat{x}_{i|k} - x^{\text{ref}}_{i|k}\rVert ^2_Q + \lVert u_{i|k}\rVert^2_R\right) + \lVert \hat{x}_{N|k} - x^{\text{ref}}_{N|k}\rVert ^2_P   \\
		 \text{s.t.} \;\;
		& \hat{x}_{i+1|k}\!=\!\!A_c(\hat{p}_{i|k})\hat{x}_{i|k}\!\!+\!\!Bu_{i|k},~i\!=\!0,\!\ldots\!,\!N\!\!-\!\!1 \\ 
		& \hat{x}_{0|k} = x_{0|k}=x_k, \\
		& \hat{x}_{i|k} \in \mathcal{X}_{i|k}, \quad \forall i = 0,1,\ldots,N,  \label{eq:NMPC_state_cons}\\
		& u_{i|k} \in \mathcal{U}_{i|k} \label{eq:NMPC_input_cons}, \quad \forall i = 0,1,\ldots,N-1.
	\end{align}
\end{subequations}
\end{problem}

The decision variables in the quadratic program (QP) in \cref{prob:Prob_LPVMPC} can be considered explicitly the control input and implicitly the states. The state and input constraints in~\cref{eq:NMPC_state_cons} and \cref{eq:NMPC_input_cons} are defined in~\cref{eq:constraints}.
For the optimization \cref{prob:Prob_LPVMPC} to be solved optimally, the estimated scheduling signal $\hat{p}_{i|k},~i=0,\ldots,N-1$, should be substituted numerically and offered prior as a prediction that will inevitably introduce error. Next, we explicitly introduce and define the error by providing the appropriate analysis for deriving deterministic bounds.
\begin{algorithm}\label{algo:lpvmpc}
\caption{The QP-based LPVMPC algorithm}
\textbf{Input}: Initial conditions $x_0$, the reference $(x^{\texttt{ref}},y^{\texttt{ref}})$ with $k\in\mathbb{Z}_+$ and the hyper-parameters $\texttt{MaxIter}\in\mathbb{Z}_+,~\varepsilon\in\mathbb{R}_+$.\\
\textbf{Output}: The control input $u_k,~k=1,\ldots$, that drives the nonlinear system to the reference under linear constraints.
\begin{algorithmic}[1]\label{alg:LPVMPC}
\State Initialize for $k=0$ the scheduling vector $\hat{p}_{i|0}$ as
$$\hat{p}_{i|0}:=\rho\left(x_0,u_0=0\right),~i=0,\ldots,N-1$$ 
\While{$k=0,1,\ldots$} 
\State Update the state $x_{i|k}^{\text{ref}}$ 
\State Set $j=0$
\While{$j<\texttt{MaxIter}$ or $\gamma_j<\varepsilon$}
\State $j\leftarrow j+1$
\State Solve the QP in \eqref{eq:LPV_MPC} 
\begin{equation*}
\begin{aligned}    
\left[\hat{x}_{i+1|k},u_{i|k}\right]&\leftarrow\texttt{QP}^{(j)}(\hat{p}_{i|k},x_k,x_{i|k}^{\text{ref}}),~i=0,\ldots,N-1\\
\text{Update}~\hat{p}_{i|k}&:=\rho(x_{i|k},u_{i|k}),~i=0,\ldots,N\\
   \gamma_j&:=\|\hat{p}_{i|k}^{(j)}-\hat{p}_{i|k}^{(j-1)}\|_2,~i=0,\ldots,N-1
\end{aligned}
\end{equation*}
  \EndWhile 
  \State Apply $u_k=u_{0|k}$ to the system \eqref{sys:LPVLQR}
  \State Measure $x_{k+1}$
  \State Update $\hat{p}_{i|k+1}=\hat{p}_{i+1|k},~i=0,\ldots,N-1$
  \State $k\leftarrow k+1$
  \EndWhile
\end{algorithmic}
\end{algorithm}
Quantifying the (parametric) uncertainty and prediction of the scheduling signal is at the heart of this study. The error source has to do with the fact that the scheduling initialization (estimation) $\hat{p}_{i|k}$ is not necessarily the correct (true) one that the nonlinear plant will follow. In particular, by implementing the optimized control input $u_{i|k}^\ast$ to the real plant, the true response $x_{i|k}$ will deviate from the predicted one $\hat{x}_{i|k}$. Thus, we aim to provide bounds for this mismatch by introducing the error state at time $k$ and prediction $i$ (i.e., at $(i+k)^{\text{th}}$ real simulation time) between the true and predicted response that is defined as
\begin{equation}\label{eq:errorstate}
    e_{i|k}:=\underbrace{x_{i|k}}_{\text{true}}-\underbrace{\hat{x}_{i|k}}_{\text{predicted}},
\end{equation}
where $e_{0|k}=x_{0|k}-\hat{x}_{0|k}=x_0-x_0=0,~\forall k\in\IZ_+$. Before deriving the theoretical analysis of the error \eqref{eq:errorstate}, we state the mean value theorem (MVT) in the generalized multivariable case, which will provide the LDI of the substantial quantities.

\begin{theorem}[The mean value theorem (MVT)]\label{the:MVT}
Let $g$ be a multivariable function defined over a nonsingleton set $\{o\}\neq[x_1,x_2]\subset\IR^n$. If
\begin{itemize}
    \item  $g$ is continuous in $[x_1,x_2]$, and
    \item  $g$ is differentiable in $(x_1,x_2)$, then
\end{itemize}
$$\exists\xi\in(x_1,x_2):g(x_2)=g(x_1)+\nabla g(\xi)(x_2-x_1).$$
\end{theorem}
\begin{remark}[Jacobian]
The Jacobian is computed
\begin{equation*}
\footnotesize
    J(x):=\nabla g(x)=\left[\begin{array}{ccc}
        \frac{\partial g_1}{\partial x_1} & \cdots &  \frac{\partial g_1}{\partial x_n} \\
        \vdots & \ddots & \vdots \\
         \frac{\partial g_n}{\partial x_1} & \cdots &  \frac{\partial g_n}{\partial x_n} \\
    \end{array}\right],~x=\left[\begin{array}{c}
         x_1  \\
         \vdots\\
         x_n
    \end{array}\right].
\end{equation*}    
\end{remark}
\section{Analysis of the error dynamics \& bounds}\label{sec:Analdyna}
At time $k$, the prediction of the scheduling signal is $\hat{p}_{i|k},~i=0,\ldots,N-1$. Schematically, we have the following high-level information:
$$
\footnotesize
\hat{p}_{i|k}\rightarrow\boxed{\hat{\Sigma}(\hat{p}_{i|k})}\xrightarrow[u^{\ast}_{i|k}]{\text{MPC}}\boxed{\Sigma(x_{i|k},u_{i|k}^\ast)}\xrightarrow[\text{to}~\Sigma]{\text{apply}} p_{i|k}=\rho(x_{i|k},u_{i|k}).
$$ 
\subsection{Error propagation within the receding horizon}
The scope of the following analysis is to provide the error propagation for the deviation between $\hat{x}_{i|k}$ (predicted state) and $x_{i|k}$ (true state) when a scheduling signal $\hat{p}_{i|k}$ has been used from a previous prediction within the receding horizon. Once a scheduling parameter has been assumed, the solution of the QP offers the optimally designed input $u^\ast$ along with the predicted states $\hat{x}$. Implementation of the input $u^\ast$ to the real system results in the true states $x$. These dynamical systems can be represented formally before and after the solution of the QP problem for $i=0,\ldots,N-1$ as  
    \begin{align}
    \footnotesize
\hat{x}_{i+1|k}&=A_c(\hat{p}_{i|k})\hat{x}_{i|k}+Bu_{i|k}^{\ast},~\text{(before QP solution)}\label{eq:LPVhat}\\
x_{i+1|k}&=A_c(p_{i|k})x_{i|k}+Bu_{i|k}^{\ast},~\text{(after QP solution)}\label{eq:LPVtrue}
    \end{align}
By having the error state from (\ref{eq:errorstate}), and after subtracting (\ref{eq:LPVhat}) from (\ref{eq:LPVtrue}), it remains a dynamical system that explains the error dynamics as
\begin{equation}\label{eq:error}
    e_{i+1|k}=A_c(p_{i|k})x_{i|k}-A_c(\hat{p}_{i|k} )\hat{x}_{i|k}.
\end{equation}
\begin{remark}[True and predicted scheduling]
    The main difference in (\ref{eq:LPVtrue}) and (\ref{eq:LPVhat}) is how the states explain the deduced scheduling signals. In particular, for the true state $x_{i|k}$, it holds $p_{i|k}=\sigma(x_{i|k},u_{i|k}^\ast)$, but for the predicted state, we cannot claim the same and the mismatch $\hat{p}_{i|k}\overset{?}{\sim}\sigma(\hat{x}_{i|k},u_{i|k}^\ast)$ produces nonzero error dynamics. 
\end{remark}
We continue with (\ref{eq:error}) after substituting the affine operator (\ref{eq:affine}) and we result to
\begin{equation}
\footnotesize
    \begin{aligned}
    e_{i+1|k}&=\left(A_{c0}+\sum_{l=1}^{n_p}p_{i|k}^{[l]}A_{cl}\right)x_{i|k}-\left(A_{c0}+\sum_{l=1}^{n_p}\hat{p}_{i|k}^{[l]}A_{cl}\right)\hat{x}_{i|k}\\
        &=A_{c0}e_{i|k}+\underbrace{\left(\sum_{l=1}^{n_p}p_{i|k}^{[l]}A_{cl}\right)}_{\sigma(x_{i|k})}x_{i|k}-\underbrace{\left(\sum_{l=1}^{n_p}\hat{p}_{i|k}^{[l]}A_{cl}\right)}_{\hat{\sigma}_{i|k}}\hat{x}_{i|k}\\
        &=A_{c0}e_{i|k}+\sigma(x_{i|k})x_{i|k}-\hat{\sigma}_{i|k}\hat{x}_{i|k},
    \end{aligned}
\end{equation}
where the smooth operator $g:\IR^n\rightarrow\IR^n$ is introduced to model the nonlinear functionality as
\begin{equation}\label{eq:operator}
    g(x):=\sigma(x)x,~x\in\IR^n,
\end{equation}
with $\sigma(x):=\sum_{l=1}^{n_p}p_{i|k}^{[l]}A_{cl}$. By denoting further the known quantity $\hat{\sigma}_{i|k}=\sum_{l=1}^{n_p}\hat{p}_{i|k}^{[l]}A_{cl}$, we can conclude to the following error dynamics $\forall k\in\IZ_+$ and $i=0,\ldots,N-1$ in a concise way as
\begin{equation}\label{eq:errordynamics}
    e_{i+1|k}=A_{c0}e_{i|k}+g(x_{i|k})-\hat{\sigma}_{i|k}\hat{x}_{i|k},~e_{0|k}:=0.
\end{equation}
The scheduling signal's initialization at time step $k$ is explained before we proceed with handling the unknown quantities in \eqref{eq:errordynamics} (i.e., $\xi_{i|k},~x_{i|k},~\hat{x}_{i|k}$).
The idea is to initialize at time $k$ with the previous true (simulated) response of the actual plant (equivalent LPV) at time $k-1$ by satisfying $\hat{p}_{i|k}=p_{i+1|k-1}=\rho(x_{i+1|k-1})=\hat{\sigma}_{i|k}$. This infers information from the previous true response of the system and provides a good prediction due to the assumed smoothness of the nonlinear operators. 
\subsection{Linear differential inclusion (LDI) \& error polytopes}
The idea is to write the LDI (exact linearization) with the help of the \cref{the:MVT} over the line segments $[x_{i+1|k-1},x_{i|k}],~i=0,\ldots,N-1,~\forall k\in\IZ_+$, which involves the two true sequential state points how to propagate over the time $k$. Applying \cref{the:MVT}, exists $\xi_{i|k}\in(x_{i+1|k-1},x_{i|k})$ such that $g(x_{i|k})=g(x_{i+1|k-1})+\nabla g(\xi_{i|k})(x_{i|k}-x_{i+1|k-1})$. Substituting in \eqref{eq:errordynamics}, it remains with $\hat{\sigma}_{i|k}=\sigma(x_{i+1|k-1})$
\begin{equation}\label{eq:errorpropagation}
\begin{aligned}
   e_{i+1|k}&=A_{c0}e_{i|k}+g(x_{i+1|k-1})-\hat{\sigma}_{i|k}\hat{x}_{i|k}\\
   &+\nabla g(\xi_{i|k})(x_{i|k}-x_{i+1|k-1}),\\
   &=A_{c0}e_{i|k}+\underbrace{\sigma(x_{i+1|k-1})}_{\hat{\sigma}_{i|k}}x_{i+1|k-1}-\hat{\sigma}_{i|k}\hat{x}_{i|k}\\
    &+\nabla g(\xi_{i|k})(x_{i|k}-x_{i+1|k-1})\Leftrightarrow\\
    e_{i+1|k}&=A_{c0}e_{i|k}+\hat{\sigma}_{i|k}(x_{i+1|k-1}-\hat{x}_{i|k})+\\
    &+\nabla g(\xi_{i|k})(x_{i|k}-x_{i+1|k-1}),~e_{0|k}=0.
    \end{aligned}
\end{equation}
Due to the uncertainty introduced by the unknowns (i.e., $\xi_{i|k},~x_{i|k},~\hat{x}_{i|k}$) in \eqref{eq:errorpropagation}, the way to proceed is to represent the error dynamics as a sequence of polytopic convex sets that will enclose the extreme behavior between the predictive and true state response deterministically due to no other disturbances. We know from the limitations of the original plant how to handle the variation of $\lVert x_{i|k}-x_{i+1|k-1}\rVert$ (e.g., physical constraints on acceleration, speed, and position). We can further enforce these constraints with the help of the decision variables $\hat{x}$ in the LPVMPC so as to impose the same variation in $\lVert x_{i+1|k-1}-\hat{x}_{i|k}\rVert$. Exploiting the affine structure of the gradient w.r.t. the scheduling parameters evaluated at $\xi_{i|k}\in(\xi_{i+1|k-1},x_{i|k})$, we can represent $\nabla g(\xi_{i|k})(x_{i|k}-x_{i+1|k-1})$ as the convex polytope $\mathbb{W}_{i|k}=Co\{\nu_{1},\nu_{2},...,\nu_{m}\},~m\in\IZ_+,~m\leq 2^{n_{\mathrm{x}}(n_{\mathrm{x}}+1)}$ where $\nu_{(\cdot)}$ represents the vertices. Similarly for the product $\hat{\sigma}_{i|k}(x_{i+1|k-1}-\hat{x}_{i|k})$, we can define the convex polytope $\mathbb{V}_{i|k}=Co\{\bar{\nu}_{1},\bar{\nu}_{2},...,\bar{\nu}_{\bar{m}}\},~\bar{m}\in\IZ_+,~\bar{m}\leq 2^{n_{\mathrm{x}}}$. By denoting the vertex $E_{0|k}:=e_{0|k}=0_{n_\textrm{x}}$, the polytopic error tubes can be computed recursively with the Minkowski summation ``$\oplus$"\footnote{Minkowski sum: $A\oplus B:=\{a+b~|~a\in A,~b\in B\}$.} and $\forall k\in\IZ_+\cup\{0\}$ along with $i=0,\ldots,N-1$ as
\begin{equation}\label{eq:errorpolytubes}
\begin{aligned}
\mathbb{E}_{i+1|k}&=\left(A_{c0}\mathbb{E}_{i|k}\right)\oplus\mathbb{V}_{i|k}\oplus\mathbb{W}_{i|k},~\mathbb{E}_{0|k}=0,\\
e_{i|k}&\in\mathbb{E}_{i|k} \Leftrightarrow x_{i|k}\in(\hat{x}_{i|k}\oplus\mathbb{E}_{i|k}).
\end{aligned}
\end{equation}
\section{Results}\label{sec:results}
\begin{example}[The unbalanced disk regulator problem]
We start by providing the \cref{tab:MPCPar} with all the input-output constraints along with the tuning parameters for solving the LPVMPC problem.
\begin{table}[h]
    \centering
    \caption{MPC Parameters}
    \label{tab:MPCPar}
   \setlength{\tabcolsep}{2pt}
    \begin{tabular}{l l | l l} 
        \textbf{Parameter}   & \textbf{Value}  & \textbf{Parameter}    & \textbf{Value}\\ 
                       Lower bound on  $\theta_k$ &   $-2\pi$ [rad] & Upper bound on  $\theta_k$ &   $2\pi$ [rad]  \\
        Lower bound on  $\omega_k$ &   $-10\pi$ [rad/s] & Upper bound on  $\omega_k$ &   $10\pi$ [rad/s] \\
        Lower bound on  $u_k$ &   $-10$ [V] & Upper bound on  $u_k$ &  $10$ [V]  \\
           Sampling time $t_s$ & $0.01$ [s] &   Horizon length $N$ & $10$    \\
          Lower bound on  $\Delta_1$ & $-t_s 10\pi$ [rad] &  Upper bound on  $\Delta_1$ & $t_s 10\pi$ [rad]\\
           Lower bound on  $\Delta_2$ & any [rad/s$^2$] &  Upper bound on  $\Delta_2$ & any [rad/s$^2$]\\
           Quadratic state costs & $Q=\texttt{diag}(8,~0.1)$ & Quadratic input cost & $R=0.5$\\
    \end{tabular}
\end{table}
The terminal cost $P$ is computed from the solution of the Lyapunov matrix equation in the robust model-based LQR together with the feedback gain $K$. The unbalanced disk control regulator problem is under consideration that starts from the initial conditions $x_0=\left[\begin{array}{cc}
    -6 & 0 \\
\end{array}\right]^\top$. We want to drive the dynamics to the state origin that is our reference signal $x^{\text{ref}}=\left[\begin{array}{cc}
   0  & 0 \\
\end{array}\right]^\top$. The dynamical system that describes the phenomenon with angular displacement-$\theta(t)$ and angular speed-$\omega(t)$ (i.e., $\omega(t)=\dot{\theta}(t)$) with state vector $x(t)=[\begin{array}{cc}
   \theta(t)  & \omega(t)   
\end{array}]^\top$ has the continuous space-state representation as in \eqref{sys:LPV} after introducing the scheduling variable $p(t):=\sin(\theta(t))/\theta(t):=\sinc(\theta(t))$. The matrices that define the continuous in-time system with $p(t)=\sinc(\theta(t))$ are
\begin{equation}
\footnotesize
    A_{cont}(p(t))=\left[\begin{array}{cc}
        0 & 1 \\
       \frac{mgl}{I_n}p(t)  & -\frac{1}{\tau}
    \end{array}\right],~B_{cont}=\left[\begin{array}{c}
         0  \\
          \frac{K_m}{\tau}
    \end{array}\right],
\end{equation}
and the chosen parameters are provided in \cref{tab:disk_parameters}.
\begin{table}
    \centering
     \caption{Parameters for the unbalanced disk example}
    \begin{tabular}{c|l}
        $I_n$ &  $2.4\cdot1e-4$ [kg$\cdot m^2$]\\
        $m$ &  $0.076~$ [Kg]\\
        $g$ & $9.81$ [m/s]\\
        $l$ & $0.041$ [m]\\
        $\tau$ & $0.4$ [1/s] \\
        $K_m$ & $11$ [rad/Vs$^2$]
    \end{tabular}
    \label{tab:disk_parameters}
\end{table}
We discretize with forward Euler\footnote{Forward Euler: $\dot{x}(t_k)\approx\frac{x(t_k+t_s)-x(t_k)}{t_s},~t_k=t_s\cdot k,~k\in\IZ_+$.}, and with $x_k=\left[\begin{matrix}
   \theta_k & \omega_k \end{matrix}\right]^\top$, the remaining dicrete LPV system \eqref{sys:LPV} with scheduling parameter $p_k=\sinc(\theta_k)$ has the following matrices
\begin{equation}
\footnotesize
\begin{aligned}\label{eq:discdisc}
    A(p_k)&=\left[\begin{array}{cc}
        1 & t_s \\
       t_s\frac{mgl}{I_n}p_k  & 1-\frac{t_s}{\tau}
    \end{array}\right],~B=\left[\begin{array}{c}
         0  \\
          t_s\frac{K_m}{\tau}
    \end{array}\right].
    \end{aligned}
\end{equation}
In the Appendix, we obtain the heavy computations of the important quantities to be used next, and we proceed from \eqref{eq:appenlast}, by defining the vector $\xi=\left[\begin{array}{cc}
   \xi^{(1)}  &  \xi^{(2)}\\
\end{array}\right]^\top$, and assume $\xi^{(1)}\in[\pi=\xi_{min}^{(1)},\xi_{max}^{(1)}=2\pi]$ where gives the maximum variation of the $-1\leq\cos(\xi^{1})\leq 1$. Then, we define
\begin{equation}
\footnotesize
    \begin{aligned}
    G_1(\xi)&=\left[\begin{array}{cc}
    0 & 0 \\
    \gamma\cos(\xi^{(1)}) & 0\end{array}\right],~h_1=\left[\begin{array}{c}
         \pm\Delta_1  \\
         \pm\Delta_2\end{array}\right],\\
         \nu_{1}^{max}&=G_1(\xi_{max}^{(1)})h_1=\left[\begin{array}{c}
              0  \\
        \gamma\cos(\xi_{max}^{(1)})\Delta_1
         \end{array}\right],\\
         \nu_{1}^{min}&=G_1(\xi_{min}^{(1)})h_1=\left[\begin{array}{c}
              0  \\
        \gamma\cos(\xi_{min}^{(1)})\Delta_1
         \end{array}\right].
    \end{aligned}
\end{equation}
Due to the switching sign in $\cos(\xi^{(1)})=\pm 1$, we can use the $+\Delta_1=t_s\omega_{max}$ only, otherwise will have the same (redundant) vertices. Moreover, it is algebraically evident that $\Delta_2$ can be any arbitrary value as it will not affect the topology of the polytope. The convex hull is defined from the two vertices only as $$G_{1}(\xi_{i|k})h_1\in\texttt{Co}(\nu_1^{min},\nu_1^{max}):=\mathbb{W}.$$
Instead of considering the most conservative case in bounding the maximum variation of the gradient, we could obtain convex sets $\mathbb{W}_{i|k}$ that can change online and take under consideration the curvature, which could reduce conservatism but it will increase the complexity. Accordingly, we define
\begin{equation}
    \footnotesize
    \begin{aligned}
        G_2(\hat{p}_{i|k})&=\left[\begin{array}{cc}
    0 & 0 \\
    \gamma\hat{p}_{i|k} & 0\end{array}\right],~h_2(\pm)=\left[\begin{array}{c}
         \pm\Delta_1  \\
         \pm\Delta_2\end{array}\right],\\
         \nu_{2}^{max}&=G_2(\hat{p}_{i|k})h_2(+)=\left[\begin{array}{c}
              0  \\
        \gamma\hat{p}_{i|k}\Delta_1
         \end{array}\right],\\
         \nu_{2}^{min}&=G_1(\hat{p}_{i|k})h_2(-)=\left[\begin{array}{c}
              0  \\
        -\gamma\hat{p}_{i|k}\Delta_1
         \end{array}\right].
    \end{aligned}
\end{equation}
Thus, the convex hull is defined from the two vertices as $$G_{2}(\hat{p}_{i|k})h_2\in\texttt{Co}(\nu_2^{min},\nu_2^{max}):=\mathbb{V}_{i|k}.$$
We denote further the vertex $\mathbb{E}_{0|k}:=e_0=\left[\begin{array}{cc}
    0 & 0 \\
\end{array}\right]^T$, and finally, the translated convex error sets around the predictions $\hat{x}_{i|k}$ with the notation in \eqref{eq:Ac0Ac1} are computed as
\begin{equation}
\begin{aligned}
    \mathbb{E}_{i+1|k}&=(A_{c0}\mathbb{E}_{i|k})\oplus\mathbb{V}_{i|k}\oplus\mathbb{W},~i=0,\ldots,N-1\\
    x_{i|k}&\in\left(\hat{x}_{i|k}\oplus\mathbb{E}_{i|k}\right),~\forall k\in\IZ_+.
\end{aligned}
\end{equation}

In \cref{fig:fig1}-(upper) is depicted the solution of the LPVMPC problem for $k=1$ along with the sequence of the error polytopes $\mathbb{E}_{i|k},~i=0,\ldots,N-1$ that gives $N$ steps ahead bounded prediction of the actual response $x_{i|k}$ compared to the predicted one. In simple terms, as the actual system has been embedded equivalently to the LPV and there are no other disturbances, we have prior the maximum variation of the true response of the system. For $k=1$, to have an improved prediction of the scheduling signal, we solve the MPC problem in \cref{alg:LPVMPC} with $\texttt{MaxIter}=10$ and $\varepsilon=1e-7$. This reduces the efficiency in the beginning as we solve a sequence of MPC problems for the fixed $k=1$, but also allows a ``warm start" that makes the initialization of the scheduling more meaningful for the system. For the rest of the simulations $k>1$, we set $\texttt{MaxIter}=1$, and we benefit from the maximum efficiency in solving QPs where the average time in (s) is $\sim 0.002<t_s=0.01$, and that certifies that we can reach real-time. 

In \cref{fig:fig1}-(lower), similarly, the aforementioned analysis is illustrated at time $k=15$.


\begin{figure*}[t]
  \centering
  \includegraphics[width=0.75\textwidth]{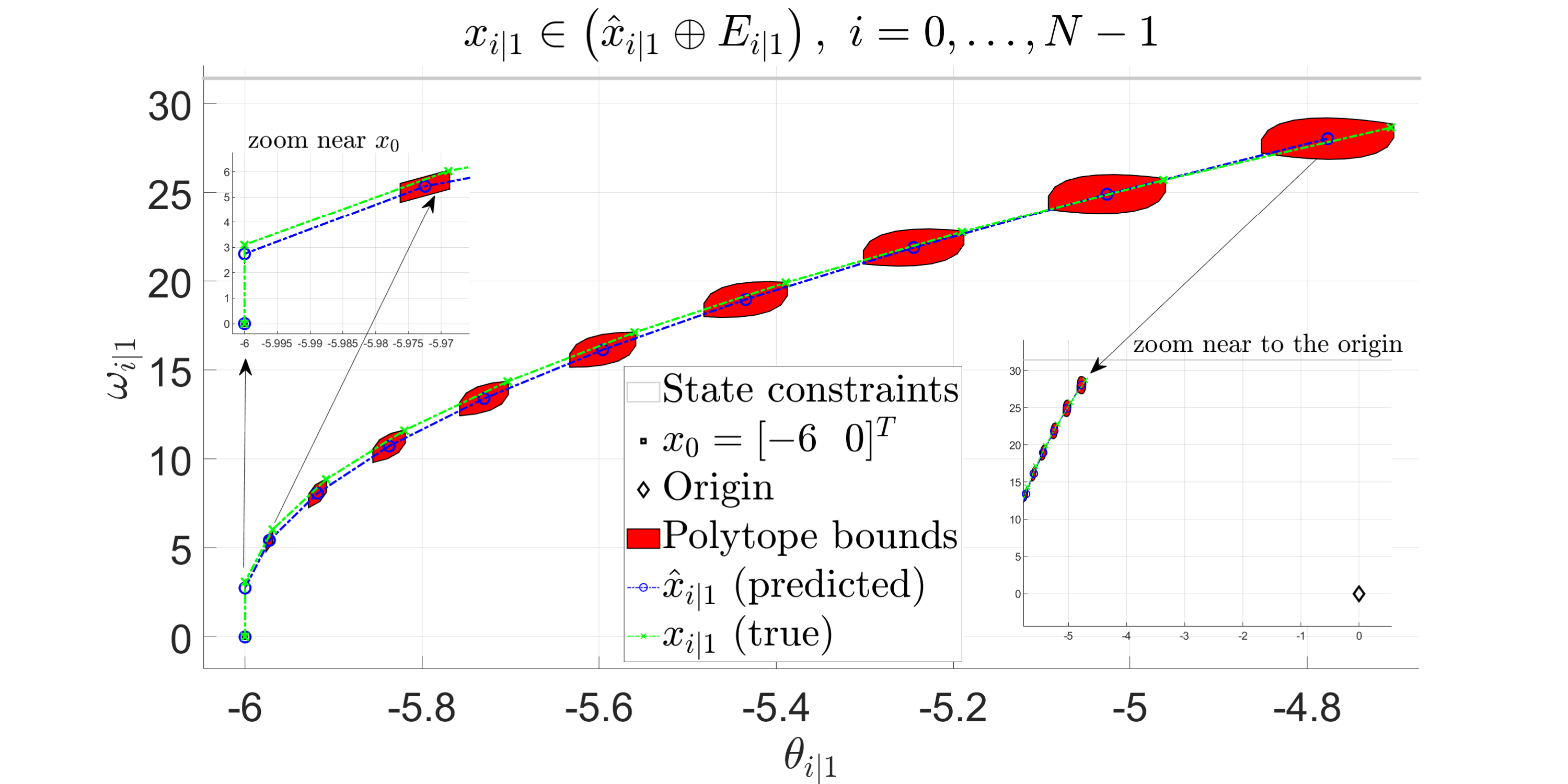}
  \hfill
  \includegraphics[width=0.75\textwidth]{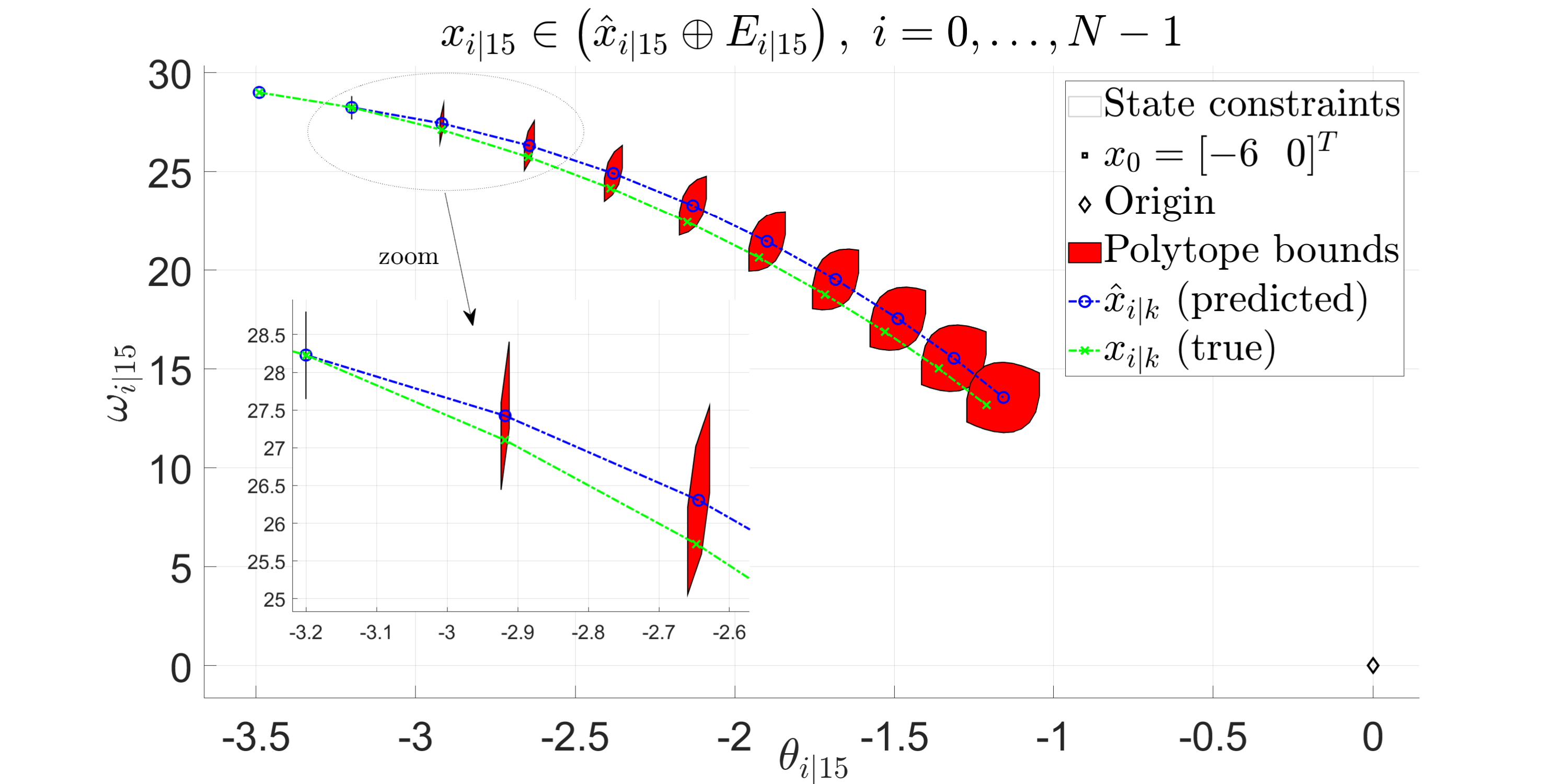}
     \caption{Upper: Phase state space evolution of the regulation problem starting from $x_0=[-6,0]^\top$ with the true (green) and predicted (blue) responses of the system at time $k=1$ and for the horizon $i=0,\ldots,10$. Error polytopes (red) centered at predicted values that enclose the true response. Lower: Phase state space predicted trajectories at $k=15$. The regulator problem solved with a fair approximation at $k=60$ which translates to real-time $0.6$ (s).}
    \label{fig:fig1}
\end{figure*}

In \cref{fig:fig3}, the complete solution with the LPVMPC framework for the regulator problem of the unbalanced disk is illustrated. The monotonicity of the angular displacement $\theta(t)$ that reaches the origin target without overshooting outlines the good performance that can be seen in NMPC frameworks. In addition, the computational burden has been avoided after utilizing the QP performance.   
\begin{figure*}[t]
    \centering
    \includegraphics[width=0.75\textwidth]{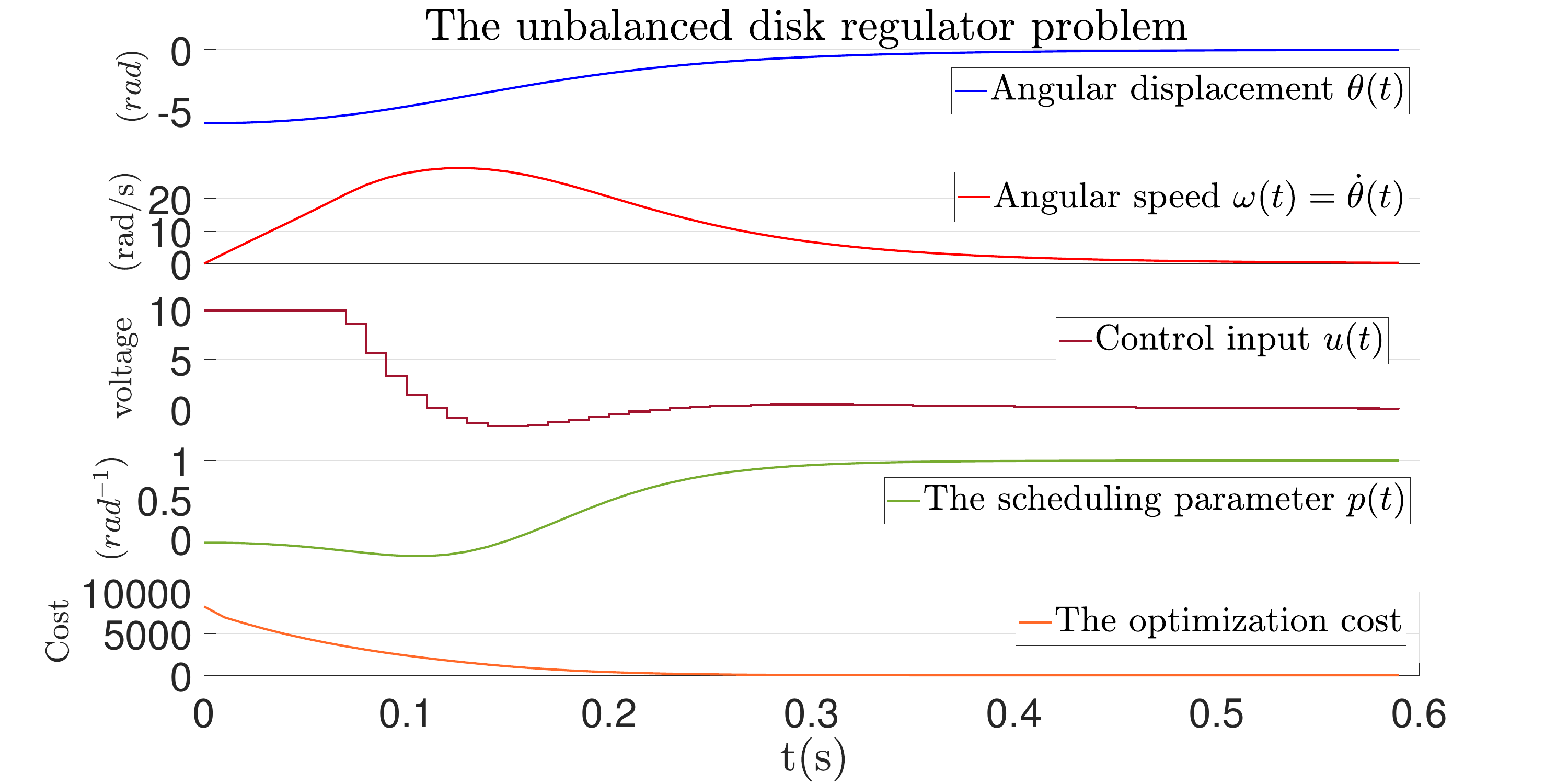}
    \caption{State, control, and scheduling trajectories of the regulation problem along with the optimization cost.}
    \label{fig:fig3}
\end{figure*}
\section{CONCLUSIONS}\label{sec:conclusion}
In this study, we were concerned with nonlinear control problems and aimed to provide a method that offers efficient solutions by analyzing some of the inherent challenges. To make MPC feasible for the nonlinear case, we represented the nonlinear system equivalently with the linear parameter-varying (LPV) embedding that can benefit from the well-established linear control theory and maintain some of the advantages like quadratic program (QP) performance robust stabilization with the linear quadratic regulator (LQR). 

The main challenge of using the LPV formulation for such control tasks was the uncertainty introduced by the predicted scheduling parameter acting as the source that dissipates errors within the receding horizon between the actual and predicted responses of the system. To tackle this challenge, we explicitly derived the error dynamics, and by applying linear differential inclusions (LDIs), we obtained the polytopic error bounds. Such a result can be computed before the MPC solution as long the actual systems' maximum capabilities (i.e., constraints) are accessible, which will help further analysis.

The derived polytopic error tubes reasonably estimate the actual response in the unbalanced disk example without being conservative. Despite that preliminary good result, in the future, we will work on reducing conservativeness further by tightening online the derived bounds by appropriately handling the polytopic sets $\mathbb{W}_{i|k},~\mathbb{V}_{i|k}$ online. In addition, we plan to use the derived error formulations with machine learning (ML) techniques, such as Gaussian processes (GPs), that can improve the predictions (e.g., to the scheduling parameters) by offering variance measures that will certify robustness against system's disturbances and measurement noise. 

We also aim to investigate more challenging systems regarding dimensionality in the state or scheduling parameter dimension. Reduction techniques will be mandatory for handling the complexity of the derived analysis. The computation of invariant sets from the derived error analysis will set the ground for our research endeavors in the immediate future. Finally, our long-term goal is to provide theoretical guarantees such as stability and recursive feasibility that will assert safety in autonomous systems that vary from mechanical to medical engineering disciplines.




\section*{APPENDIX}
Stabilizing the unstable plant with the robust model-based LQR feedback, the gain $K$ is computed optimally. To simplify the exposition, suppose that $K=\left[\begin{array}{cc}
    \alpha & \beta\end{array}\right]$, thus the stabilized linear matrix $A_c(p_k)$ has the following form:
    \begin{equation*}
    \footnotesize
    \begin{aligned}
     A_c(p_k)&=\left[\begin{array}{cc}
        1 & t_s \\
        t_s\frac{mgl}{I_n}p_k+\alpha t_s\frac{K_m}{\tau} & 1-\frac{t_s}{\tau}+\beta t_s\frac{K_m}{\tau} 
    \end{array}\right]
    \end{aligned}
    \end{equation*}
    We define the parameters $\gamma:=t_smgl/I_n$, $\delta:=\alpha t_sK_m/\tau$, $\eta:=1-t_s/\tau+\beta t_sK_m/\tau$, and exploit the affine structure
    \begin{equation}\label{eq:Ac0Ac1}
    \footnotesize
        A_c(\rho(\theta_k))=\underbrace{\left[\begin{array}{cc}
            1 & t_s \\
            \delta & \eta 
        \end{array}\right]}_{A_{c0}}+\underbrace{\left[\begin{array}{cc}
            0 & 0 \\
            \gamma  & 0
        \end{array}\right]}_{A_{c1}}\sinc(\theta_k).
    \end{equation}
    The operator $g$ with $x=\left[\begin{array}{cc}
         \theta & \omega \\
     \end{array}\right]^T$  remains
    \begin{equation*}
\footnotesize
    \begin{aligned}
        g(x)&=\sigma\left(\left[\begin{array}{c}
             \theta  \\
             \omega 
        \end{array}\right]\right)\left[\begin{array}{c}
         \theta  \\
         \omega
    \end{array}\right]=\sinc(\theta)\left[\begin{array}{c}
         \theta  \\
         \omega
    \end{array}\right]=\left[\begin{array}{c}
         \sin(\theta)  \\
         \omega\sinc(\theta) 
    \end{array}\right],
    \end{aligned}
    \end{equation*}
    where the Jacobian is computed as
    \begin{equation*}
    \footnotesize
        J(x)=\nabla g(x)=\left[\begin{array}{cc}
        \cos(\theta) & 0\\
        \frac{\omega}{\theta}\left(\cos(\theta)-\sinc(\theta)\right) & \sinc(\theta)
    \end{array}\right].
    \end{equation*}
\end{example}
Substituting in \eqref{eq:errorpropagation}, we can derive explicitly
\begin{equation}\label{eq:appenlast}
\footnotesize
    \begin{aligned}
        e_{i+1|k}&=\left[\begin{array}{cc}
            1 & t_s \\
            \delta & \eta
        \end{array}\right]e_{i|k}+\left[\begin{array}{cc}
           0 & 0 \\
           \gamma\hat{p}_{i|k}  & 0
        \end{array}\right]\left[\begin{array}{c}
             x^{(1)}_{i+1|k-1}-\hat{x}^{(1)}_{i|k}  \\
             x^{(2)}_{i+1|k-1}-\hat{x}^{(2)}_{i|k} 
        \end{array}\right]\\
        &+\left[\begin{array}{cc}
           0 & 0 \\
           \gamma\cos(\xi_{i|k}^{(1)})  & 0
        \end{array}\right]\left[\begin{array}{c}
             x^{(1)}_{i|k}-\hat{x}^{(1)}_{i+1|k-1}  \\
             x^{(2)}_{i|k}-\hat{x}^{(2)}_{i+1|k-1}\end{array}\right].
    \end{aligned}
\end{equation}
%
%
\bibliographystyle{IEEEtran}
\bibliography{main}
\end{document}